%% file: AMSReesDet.tex
\documentclass[12pt,reqno]{amsart}
\usepackage{amssymb,latexsym, mathtools,tikz}
\usepackage{amsfonts}
\usepackage{amsmath}
\usepackage{amsthm}
\usepackage{amsxtra}
\usepackage{bbm}
\usepackage{braket}
\usepackage{comment}
\usepackage{extarrows}
\usepackage{graphics}
\usepackage{latexsym}
\usepackage{float}
\usepackage[colorlinks=true,citecolor=blue,linkcolor=magenta]{hyperref} 
\usepackage{cleveref}

\usepackage[hyperpageref]{backref}
\usepackage{enumerate}
\usepackage[top=35mm, bottom=25mm, left=30mm, right=30mm]{geometry}

\begin{document}

\title{Rees Algebras of unit interval determinantal facet ideals}


\author{Ayah Almousa}
\address{Department of Mathematics, University of Minnesota - Twin Cities, Minneapolis, MN 55455}
\email{almou007@umn.edu}

\author{Kuei-Nuan Lin}
\address{Department of Mathematics, Penn State University, Greater Allegheny, McKeesport, PA 15132}
\email{linkn@psu.edu}

\author{Whitney Liske}
\address{Department of Mathematics, Saint Vincent College, Latrobe, PA 15650}
\email{whitney.liske@stvincent.edu}





\dedicatory{}


\begin{abstract}
Using SAGBI basis techniques, we find Gr\"obner bases for the presentation ideals of the Rees algebras and special fiber rings of unit interval determinantal facet ideals. In particular, we show that unit interval determinantal facet ideals are of fiber type and that their special fiber rings are Koszul. Moreover, their Rees algebras and special fiber rings are normal Cohen-Macaulay domains and have rational singularities. 
\end{abstract}

\maketitle

\input macro.tex

\newcommand{\oR}{\overline{R}}
\newcommand{\bfu}{{\mathbf u}}
\renewcommand{\nn}{{\mathbf n}}
\newcommand{\bfc}{{\mathbf c}}
\renewcommand{\ayah}[1]{{\color{magenta} \sf AYAH: [#1]}}
\renewcommand{\supp}{\text{Supp}}
\input intro
\input background

\input reesClosed
\bibliographystyle{amsplain}
\bibliography{bib}
\addcontentsline{toc}{section}{Bibliography}
\end{document}

%% file: macro.tex
\theoremstyle{plain}
\newtheorem{thm}{Theorem}[section]
\newtheorem{lemma}[thm]{Lemma}
\newtheorem{cor}[thm]{Corollary} 
\newtheorem{prop}[thm]{Proposition}
\newtheorem{exercise}[thm]{Exercise}
\newtheorem{claim}[thm]{Claim}
\newtheorem{conjecture}[thm]{Conjecture}
\newtheorem{question}[thm]{Question}
\newtheorem{problem}[thm]{Problem}
\newtheorem{openpr}[thm]{Open-Ended Problem}

\theoremstyle{definition}
\newtheorem{definition}[thm]{Definition}
\newtheorem{example}[thm]{Example}
\newtheorem{examples}[thm]{Examples}
\newtheorem{notation}[thm]{Notation}
\newtheorem{assumptions}[thm]{Assumptions}
\newtheorem{motivation}[thm]{Motivation}
\newtheorem{remark}[thm]{Remark}
\newtheorem{observation}[thm]{Observation}
\newtheorem{construction}[thm]{Construction}
\newtheorem{pr}[thm]{Project}
\newtheorem*{acknowledgment*}{Acknowledgments}
\def\citem#1{\item[{\rm (#1)}]} 
\def\kk{\mathbb{K}}		

\newenvironment{bmlist}    
{\begin{list}{\text} 
{\setlength{\rightmargin}{0cm}\setlength{\itemindent}{0cm}
\setlength{\labelwidth}{1cm}\setlength{\leftmargin}{1.2cm}
\setlength{\labelsep}{.2cm}}}{\end{list}} 

\def\bprop{\begin{prop}}\def\eprop{\end{prop}}
\def\bthm{\begin{thm}}\def\ethm{\end{thm}}
\def\bcor{\begin{cor}}\def\ecor{\end{cor}}
\def\blemma{\begin{lemma}}\def\elemma{\end{lemma}}
\def\bexercise{\begin{exercise}}\def\eexercise{\end{exercise}}
\def\bconjecture{\begin{conjecture}}\def\econjecture{\end{conjecture}}
\def\bclaim{\begin{claim}}\def\eclaim{\end{claim}}
\def\bquestion{\begin{question}}\def\equestion{\end{question}}
\def\bproblem{\begin{problem}}\def\eproblem{\end{problem}}
\def\boproblem{\begin{openpr}}\def\eoproblem{\end{openpr}}

\def\bconstruction{\begin{construction}}\def\econstruction{\end{construction}}
\def\bdefinition{\begin{definition}}\def\edefinition{\end{definition}}
\def\bremark{\begin{remark}}\def\eremark{\end{remark}}
\def\bnotation{\begin{notation}}\def\enotation{\end{notation}}
\def\bobservation{\begin{observation}}\def\eobservation{\end{observation}}
\def\bexample{\begin{example}}\def\eexample{\end{example}}
\def\b{\begin}   
\def\l{\label}    

\def \bspace{\vglue .2cm}
\def\mhbreak{\hfil\break}
\def\mvbreak{\vfil\break}
\def\pd{{\rm{pd}}}

\def\bfix#1{{{\bf***FIX:} #1 {\bf **}}}

\DeclarePairedDelimiter\abs{\lvert}{\rvert}%

\newcommand{\initial}{\operatorname{in}}
\newcommand{\NF}{\operatorname{NF}}
\newcommand{\HF}{\operatorname{HF}}
\newcommand{\Hilb}{\operatorname{Hilb}}
\newcommand{\depth}{\operatorname{depth}}
\newcommand{\reg}{\operatorname{reg}}
\newcommand{\Span}{\operatorname{span}}
\newcommand{\img}{\operatorname{img}}
\newcommand{\inn}{\operatorname{in}}
\newcommand{\lcm}{\operatorname{lcm}}

\newcommand{\sgn}{\operatorname{sgn}}
\newcommand{\length}{\operatorname{length}}
\newcommand{\coker}{\operatorname{coker}}
\newcommand{\adeg}{\operatorname{adeg}}
\newcommand{\pdim}{\operatorname{pdim}}
\newcommand{\Spec}{\operatorname{Spec}}
\newcommand{\Ext}{\operatorname{Ext}}
\newcommand{\Tor}{\operatorname{Tor}}
\newcommand{\LT}{\operatorname{LT}}
\newcommand{\im}{\operatorname{im}}
\newcommand{\NS}{\operatorname{NS}}
\newcommand{\Frac}{\operatorname{Frac}}
\newcommand{\Char}{\operatorname{char}}
\newcommand{\Proj}{\operatorname{Proj}}
\newcommand{\id}{\operatorname{id}}
\newcommand{\Div}{\operatorname{Div}}
\newcommand{\Cl}{\operatorname{Cl}}
\newcommand{\tr}{\operatorname{tr}}
\newcommand{\Tr}{\operatorname{Tr}}
\newcommand{\ann}{\operatorname{ann}}
\newcommand{\Gal}{\operatorname{Gal}}
\newcommand{\QQbar}{{\overline{\mathbb Q}}}
\newcommand{\Br}{\operatorname{Br}}
\newcommand{\Bl}{\operatorname{Bl}}
\newcommand{\Cox}{\operatorname{Cox}}
\newcommand{\conv}{\operatorname{conv}}
\newcommand{\getsr}{\operatorname{Tor}}
\newcommand{\diam}{\operatorname{diam}}
\newcommand{\Hom}{\operatorname{Hom}} 
\newcommand{\sheafHom}{\mathcal{H}om}
\newcommand{\Gr}{\operatorname{Gr}}
\newcommand{\rank}{\operatorname{rank}}
\newcommand{\codim}{\operatorname{codim}}
\newcommand{\Sym}{\operatorname{Sym}} 
\newcommand{\GL}{{GL}}
\newcommand{\Prob}{\operatorname{Prob}}
\newcommand{\Density}{\operatorname{Density}}
\newcommand{\Syz}{\operatorname{Syz}}
\newcommand{\Jac}{\operatorname{Jac}}
\renewcommand{\pd}{\operatorname{pd}}
\newcommand{\supp}{\operatorname{supp}}
\newcommand{\cone}{\operatorname{\textbf{cone}}}
\newcommand{\Res}{\operatorname{Res}}
\newcommand{\deghat}{\widehat{\deg}}
\newcommand{\red}{\operatorname{red}}
\newcommand{\ord}{\operatorname{ord}}
\newcommand{\sort}{\operatorname{sd}}
\newcommand{\inv}{\operatorname{inv}}

\newcommand{\defi}[1]{\textsf{#1}} 

\newcommand{\colim}{\operatorname{colim}}
\newcommand{\doot}{\bullet}

\newcommand{\Alt}{\bigwedge\nolimits}
\newcommand{\Sch}{\text{\bf Sch}}										
\newcommand{\clique}{\Delta^{\textrm{clique}}}


\renewcommand{\aa}{\mathbf a}
\newcommand{\bb}{\mathbf b}
\newcommand{\cc}{\mathbf c}
\newcommand{\dd}{\mathbf d}
\newcommand{\ee}{\mathbf e}
\newcommand{\ff}{\mathbf f}
\renewcommand{\gg}{\mathbf g}

\newcommand{\mm}{\mathbf m}
\newcommand{\nn}{\mathbf n}
\newcommand{\pp}{\mathbf p}
\newcommand{\xx}{\mathbf x}
\newcommand{\uu}{\mathbf u}
\newcommand{\vv}{\mathbf v}
\newcommand{\ww}{\mathbf w}
\newcommand{\yy}{\mathbf y}

\newcommand{\bfF}{\mathbf F}
\newcommand{\bfG}{\mathbf G}
\newcommand{\bfE}{\mathbf E}
\newcommand{\bfM}{\mathbf M}

\newcommand{\bfN}{\mathbf N}
\newcommand{\bfQ}{\mathbf Q}

\newcommand{\hH}{\operatorname{h}}
\renewcommand{\H}{\operatorname{H}}
\newcommand{\K}{\operatorname{K}}
\newcommand{\kK}{\operatorname{k}}
\newcommand{\OO}{\operatorname{O}}
\newcommand{\oo}{\operatorname{o}}

\newcommand{\cA}{\mathcal{A}}
\newcommand{\cB}{\mathcal{B}}
\newcommand{\cC}{\mathcal{C}}
\newcommand{\cE}{\mathcal{E}}
\newcommand{\cF}{\mathcal{F}}
\newcommand{\cG}{\mathcal{G}}
\newcommand{\cH}{\mathcal{H}} 
\newcommand{\cI}{\mathcal{I}}
\newcommand{\cJ}{\mathcal{J}}
\newcommand{\cK}{\mathcal{K}}
\newcommand{\cL}{\mathcal{L}}
\newcommand{\cM}{\mathcal{M}}
\newcommand{\cN}{\mathcal{N}}
\renewcommand{\O}{\mathcal{O}}
\newcommand{\cP}{\mathcal{P}}
\newcommand{\cQ}{\mathcal{Q}}
\newcommand{\cR}{\mathcal{R}}
\newcommand{\cS}{\mathcal{S}}
\newcommand{\cT}{\mathcal{T}}
\newcommand{\U}{\mathcal{U}} 		
\newcommand{\cV}{\mathcal{V}}
\newcommand{\cW}{\mathcal{W}}
\newcommand{\cX}{\mathcal{X}}
\newcommand{\cY}{\mathcal{Y}}
\newcommand{\cZ}{\mathcal{Z}}

\newcommand{\BB}{\mathbb{B}}
\newcommand{\CC}{\mathbb{C}}
\newcommand{\DD}{\mathbb{D}}
\newcommand{\EE}{\mathbb{E}}
\newcommand{\FF}{\mathbb{F}}
\newcommand{\GG}{\mathbb{G}}
\newcommand{\II}{\mathbb{I}}
\newcommand{\JJ}{\mathbb{J}}
\newcommand{\LL}{\mathbb{L}}
\newcommand{\MM}{\mathbb{M}}
\newcommand{\NN}{\mathbb{N}}
\newcommand{\PP}{\mathbb{P}}
\newcommand{\QQ}{\mathbb{Q}}
\newcommand{\RR}{\mathbb{R}}
\renewcommand{\S}{\mathbb{S}}
\newcommand{\TT}{\mathbf{T}}
\newcommand{\bF}{\mathbf{F}}
\newcommand{\UU}{\mathbb{U}}		
\newcommand{\VV}{\mathbb{V}}
\newcommand{\WW}{\mathbb{W}}
\newcommand{\XX}{\mathbb{X}}
\newcommand{\YY}{\mathbb{Y}}
\newcommand{\ZZ}{\mathbb{Z}}

\newcommand{\sA}{\mathsf{A}}
\newcommand{\sB}{\mathsf{B}}
\newcommand{\sC}{\mathsf{C}}
\newcommand{\sD}{\mathsf{D}}
\newcommand{\sE}{\mathsf{E}}
\newcommand{\sF}{\mathsf{F}}
\newcommand{\sG}{\mathsf{G}}
\newcommand{\sH}{\mathsf{H}} 
\newcommand{\sI}{\mathsf{I}}
\newcommand{\sJ}{\mathsf{J}}
\newcommand{\sK}{\mathsf{K}}
\newcommand{\sL}{\mathsf{L}}
\newcommand{\sM}{\mathsf{M}}
\newcommand{\sN}{\mathsf{N}}
\newcommand{\sO}{\mathsf{O}}
\newcommand{\sP}{\mathsf{P}}
\newcommand{\sQ}{\mathsf{Q}}
\newcommand{\sR}{\mathsf{R}}
\newcommand{\sS}{\mathsf{S}}
\newcommand{\sT}{\mathsf{T}}
\newcommand{\sU}{\mathsf{U}} 
\newcommand{\sV}{\mathsf{V}}
\newcommand{\sW}{\mathsf{W}}
\newcommand{\sX}{\mathsf{X}}
\newcommand{\sY}{\mathsf{Y}}
\newcommand{\sZ}{\mathsf{Z}}
 
\newcommand{\cl}{\mathfrak{cl}}
\newcommand{\g}{\mathfrak{g}}
\newcommand{\h}{\mathfrak{h}}
\newcommand{\m}{\mathfrak{m}}
\newcommand{\n}{\mathfrak{n}}
\newcommand{\p}{\mathfrak{p}}
\newcommand{\q}{\mathfrak{q}}
\renewcommand{\r}{\mathfrak{r}}

\newcommand{\ayah}[1]{{\color{magenta} \sf AYAH: [#1]}}
\newcommand{\lin}[1]{{\color{blue} \sf LIN: [#1]}}
\newcommand{\whit}[1]{{\color{orange} \sf Whitney: [#1]}}

\newcommand{\ra}{\rightarrow}

\newcommand{\Xv}{\check{X}}
\newcommand{\Yv}{\check{Y}}
\newcommand{\Zv}{\check{Z}}

%% file: intro.tex
\section{Introduction}
In this work, we study the blow-ups of certain determinantal varieties called \emph{determinantal facet ideals}. To be more specific, we find the homogeneous coordinate rings of graphs and images of the blow-ups of a projective space along its subscheme defined by a certain class of determinantal varieties.
Given an ideal $I$ in a polynomial ring $R=\mathbb{K}[x_1,...,x_n]$ over a field $\mathbb{K}$, the \emph{Rees algebra} of $I$ is defined to be the graded algebra $\mathcal{R}(I) = \oplus_{i=0}^{\infty}I^{i}t^{i}\subset R[t]$, where $t$ is an indeterminate over $R$. The \emph{special fiber ring} $\mathcal{F}(I)$ is defined as $\mathcal{R}(I)\otimes \mathbb{K}$. The projective schemes of $\mathcal{R}(I)$ and $\mathcal{F}(I)$ define the blowup and the special fiber of the blowup of the scheme $\mathrm{Spec}(R)$ along $V(I)$.

The Rees algebra is an important object in commutative algebra, algebraic geometry, elimination theory, intersection theory, geometric modeling, chemical reaction networks, and many more fields; see \cite{Cox} and \cite{Cox-Lin-Sosa} for details on such applications.
If the ideal $I$ is minimally generated by $\mu$ elements, we find ideals $\mathcal{J}$ and $\mathcal{K}$ over polynomial rings $S = R[T_1 ,\ldots, T_{\mu}]$ and $\kk[\TT] = \mathbb{K}[T_1,\ldots,T_{\mu}]$ respectively, such that $\mathcal{R} (I) = S/\mathcal{J}$ and $\mathcal{F} (I) = \kk[\TT] /\mathcal{K}$.
The defining equations of $\mathcal{J}$ and $\cK$ are implicit equations of the varieties defined by the graph and image of a blow-up, respectively. Finding the implicit equations of the presentation ideals $\cJ$ and $\cK$ of $\mathcal{R}(I)$ and $\mathcal{F}(I)$, respectively, is a challenging problem and is still open for many classes of ideals. In particular, the presentation ideals of the Rees algebra of determinantal ideals are only known in very special cases. Conca, Herzog, and Valla found the presentation ideal for the Rees algebra of the ideal of maximal minors of a generic matrix in \cite{CHV96}. The presentation ideal for the Rees algebra of the ideal of the rational normal scroll associated with a $2\times n$ matrix was shown by Sammartano in \cite{sammartano2020blowup}. Very recently, the case of two-minors of a generic $3\times n$ matrix was resolved by Huang, Perlman, Polini, Raicu, and Sammartano in \cite{notreDame2minorRels}. 

Let $X = (x_{ij})$ be a generic $m\times n$ matrix over a ring $R=\mathbb{K}[X]$, and assume $m\leq n$. It is well-known that the ideal of maximal minors of $X$, denoted $I_m(X)$, is of fiber type; that is, the ideal $\mathcal{J}$ is generated by linear relations with respect to the variables in $\TT$, and the generators of $\mathcal{K}$ (see \cite{CHV96} and Subsection \ref{subsec: rees}).
The presentation ideal of the special fiber ring, $\mathcal{K}$, is the ideal of a Grassmannian defined by Pl\"ucker relations; see, for example, \cite[Chapter 14]{MillerSturmfels}. Moreover, Eisenbud and Huneke proved in \cite{eisenbudHuneke1983} that the Rees algebra of $I_m(X)$ is a normal Cohen-Macaulay domain.
In the 1980s, Bruns, Simis, and Trung considered the ideal generated by all maximal minors which share the first $k$ columns of $X$ and showed that they are always of fiber type in \cite{brunsSimisTrung1991Hodge}. Bruns and Simis found the symmetric algebra for this class of ideals in \cite{brunsSimis1987symmetric}. In later work with Trung, they used the Hodge algebra structure on these ideals to give the defining equations of the Rees algebra in \cite{brunsSimisTrung1991Hodge} and concluded that they are of fiber type.
To the best of the authors' knowledge, not much more is known for Rees algebras of sub-ideals of $I_m(X)$. Even in the case when the ideal is generated by a subset of maximal minors of a $2\times 5$ matrix, the ideal may not be of fiber type; see Example \ref{ex: notFiberType}. Clearly, one needs to impose extra conditions in order to have hope of describing generators of presentation ideals of Rees algebras and their properties. 

Determinantal facet ideals, which were introduced by Ene, Herzog, Hibi, and Mohammadi in \cite{EHHM13}, are generated by a subset of maximal minors of an $m\times n$ matrix indexed by the facets of a pure $(m-1)$-dimensional simplicial complex $\Delta$ on $n$ vertices. They are a natural generalization of binomial edge ideals. Recall that binomial edge ideals were introduced by Herzog, et. al. in \cite{BEIpaper} due to their connections with algebraic statistics; see, for example, \cite{diaconis1998lattice}.
Particularly tractable classes of determinantal facet ideals arise from minimal generating sets which are indexed by the facets of the corresponding simplicial complex. They also form a Gr\"obner basis with respect to the lexicographic monomial order $>$ induced by $x_{11} > x_{12} > \cdots > x_{1n} > x_{21} > x_{22} > \cdots >x_{mn}$. Ideals generated by maximal minors of a matrix are well-known to form a Gr\"obner basis with respect to \textit{any} term order (see \cite{SturmfelsZelevinsky93}). The authors of \cite{BEIpaper} characterized binomial edge ideals for which the quadratic generators form a Gr\"obner basis with respect to $>$; this class of binomial edge ideals is called \textit{closed}.
Closed binomial edge ideals have been extensively studied by many authors; see \cite{saeedi2016binomial} for a compilation of results on closed binomial edge ideals.
There have been several attempts by various authors, including Almousa--Vandebogert \cite{almousa2022determinantal} and Benedetti--Seccia--Varbaro \cite{benedetti2021hamiltonian} to find a necessary and sufficient condition for the natural minimal generating set of determinantal facet ideals to form a Gr\"obner basis with respect to $>$, but it has proven to be a difficult task.

One class of determinantal facet ideals that Almousa--Vandebogert and Benedetti--Seccia--Varbaro independently settled on as a natural candidate to study is \textit{unit interval} determinantal facet ideals (see Definition \ref{def: unit}), a class that has a natural combinatorial description and generalizes closed binomial edge ideals.
The natural minimal generating set of a unit interval determinantal facet ideal corresponding to the facets of a simplicial complex is a reduced Gr\"obner basis with respect to any diagonal term order; see \cite{almousa2022determinantal}. Therefore, it is natural to study the homological properties of this class of ideals by first aiming to understand the homological properties of their initial ideals.
Ene, Herzog, and Hibi conjectured in \cite{ene2011cohen} that the graded Betti numbers of a closed binomial edge ideal and its initial ideal with respect to a diagonal term order coincide. Almousa and Vandebogert asked whether this might extend to all \textit{lcm-closed} determinantal facet ideals, a class that includes unit interval determinantal facet ideals. By studying the linear strands of initial ideals of lcm-closed determinantal facet ideals, this was confirmed in \cite{almousa2020resolutions} for the case when the corresponding simplicial complex has no more than two maximal cliques.

We aim to study the Rees algebra of a unit interval determinantal facet ideal via the study of the Rees algebra of its initial ideal using the theory of SAGBI bases.
This approach was used successfully by Conca, Herzog, and Valla in \cite{CHV96} to find the presentation ideal of the Rees algebra for balanced rational normal scrolls, and by Lin and her coauthors for the secant varieties of balanced rational normal scrolls in \cite{LinShen}, and sparse matrices in \cite{WICA}. The Cohen-Macaulayness of rational normal scrolls is shown in \cite{LinShen2}.

The paper is outlined as follows. We establish notation and recall some preliminaries in Section \ref{sec:prelims}.
In Section \ref{sec: Clique-Sorted}, we present the candidates for the defining equations of the Rees algebra of the initial ideal of a unit interval determinantal facet ideal in \Cref{reesEQ}. We define the ``clique-sorted'' monomials in \Cref{def: cliqueSorted}, and show that these clique-sorted monomials are unique modulo the ideal generated by the candidates in \Cref{Unique-CliqueSorted}. We conclude that polynomials described in \Cref{reesEQ} form a Gr\"obner basis for the Rees algebra of the initial ideal of a unit interval determinantal facet ideal in Theorem \ref{thm: ReesInitialClosedDFI}, extending a result of Ene, Herzog, Hibi, and Mohammadi in \cite[Corollary 1.4]{EHHM13}. Then we show that these polynomials given in \Cref{reesEQ} can be lifted to a Gr\"obner basis for the Rees algebra of a unit interval determinantal facet ideal in Theorem \ref{thm: reesClosedDFI}. Therefore any unit interval determinantal facet ideal is of fiber type. We further give necessary and sufficient conditions for it to be of linear type, recovering a theorem of Bruns, Simis, and Trung in \cite{brunsSimisTrung1991Hodge}. In particular, the special fiber ring of any unit interval determinantal facet ideal is Koszul and its presentation ideal is generated by Pl\"ucker relations in Corollary \ref{cor: propertiesReesClosed}. Finally, via the SAGBI basis deformation, we see that both the Rees algebra and special fiber ring of a unit interval determinantal facet ideal are normal Cohen-Macaulay domains and have rational singularities, extending a result of Eisenbud and Huneke in \cite{eisenbudHuneke1983}.

%% file: background.tex
\section{Preliminaries}\label{sec:prelims}
\subsection{Determinantal Facet Ideals}

Let $X = (x_{ij})$ be an $m\times n$ matrix of indeterminates where $m\leq n$, and let $R = 
\mathbb{K}[X]$ be the polynomial ring over a field $\mathbb{K}$ in the indeterminates $x_{ij}$. For indices $\aa = \{a_1,\ldots, a_m\}$ such that $1\leq a_1 < \cdots < a_m \leq n$, set $[\aa]$ = $[a_1 \cdots a_m]$ to be the maximal minor of $X$ involving columns in $\aa$.  The ideal generated by all $m$-minors of $X$ is denoted by $I_m(X)$. 

\begin{definition}
\begin{enumerate} [(a)]
\item \label{def: DFI} Let $\Delta$ be a pure $(m-1)$-dimensional simplicial complex on the vertex set $V=\left[n\right]$. A \emph{determinantal facet ideal} (or \textit{DFI}) $J_\Delta\subseteq R$ is the ideal generated by determinants of the form $[\aa]$ where $\aa$ supports an $m-1$ face of $\Delta$; that is, the columns of $[\aa]$ correspond to a facet $F=\{a_1,\ldots a_m\} \in \Delta$. When $m=2$, one may identify $\Delta$ with a graph $G$ and $J_G$ is called a \emph{binomial edge ideal}.
\item \label{def:cliques} Let $\Delta$ be an $(m-1)$-dimensional simplicial complex on vertex set $[n]$.
A \textit{clique} of $\Delta$ is an induced subcomplex $\Gamma$ of $\Delta$ such that any $m$ vertices of $\Gamma$ are in a face together.
A clique is called \textit{maximal} if it is not contained in any larger clique of $\Delta$.
The simplicial complex $\clique$ whose facets are the maximal cliques of $\Delta$ is called the \emph{clique complex} associated to $\Delta$. The decomposition $\Delta = \Delta_1 \cup \cdots \cup \Delta_r$ is called the \emph{maximal clique decomposition} of $\Delta$ where $\Delta_i$'s are maximal cliques of $\Delta$. 

\end{enumerate} 
\end{definition}

\begin{remark} Let $I$ be an ideal generated by an arbitrary subset of maximal minors of $X$. The simplicial complex $\Delta$ associated to a determinantal facet ideal can be viewed as a combinatorial tool to index generators of such an ideal, since the vertices of each facet correspond to the columns defining a minor in the generating set of $I$. For any $\Delta_i$ in the clique decomposition of $\Delta$, let $V_i$ denote the vertex set of $\Delta_i$. Then each $\Delta_i$ corresponds to a submatrix $X_{\Delta_i}$ of $X$ with columns in the set $V_i$ such that the ideal of maximal minors $I_m(X_{\Delta_i})$ is contained in $J_\Delta$. 
\end{remark}

\begin{notation}\label{not: initialIdeal} Let $>$ denote the lexicographic monomial order induced by the natural order of indeterminates $x_{11} > x_{12} > \cdots >x_{21} > x_{22} > \cdots > x_{mn}$. 
\begin{enumerate}[(i)]
    \item Set $\xx_\aa = \inn_>[\aa] = x_{1a_1}x_{2a_2}\cdots x_{ma_m}$. Frequently, we will drop $>$ and simply write $\inn(J_\Delta)$ for the initial ideal of $J_\Delta$.
    \item Set $\TT_\Delta = \{T_\aa \mid \aa \text{ is an $(m-1)$-face of $\Delta$}\}$.
\end{enumerate}

\end{notation}

In general, the initial ideal of a determinantal facet ideal with respect to an arbitrary term order is not well-understood outside of the case when $m=2$.
However, one class where the Gr\"obner basis of a determinantal facet ideal with respect to $>$ is well-understood is a \textit{unit interval} DFI. In this case, the corresponding simplicial complex is a \textit{unit interval} simplicial complex. These were introduced independently in \cite{almousa2022determinantal} and \cite{benedetti2021hamiltonian}.

\begin{definition}\label{def: unit} Let $\Delta$ be a pure $(m-1)$-dimensional simplicial complex on $n$ vertices with maximal clique decomposition $\Delta = \bigcup_{i=1}^r \Delta_i$. 
The simplicial complex $\Delta$ is a \emph{unit interval} simplicial complex if each $\Delta_i$ may be written as an interval $[u_i , v_i] = \{ u_i , u_i+1 , \dots , v_i-1 , v_i \}$ for integers $u_i < v_i$.  We may assume $$1=u_1< u_2 < \cdots < u_r<v_n=n.$$
We call the determinantal facet ideal $J_\Delta$ unit interval if the corresponding simplicial complex $\Delta$ is a unit interval simplicial complex.
\end{definition}

\begin{figure}[h]
\begin{tikzpicture}[scale=0.4]
\filldraw [gray!50](-1,0)--(4,4)--(6,-1)--(3,-3)--(-1,0);
\filldraw [gray!80](4,4)--(3,-3)--(2,1)--(4,4);
\draw [shading=ball, ball color=black]  (-1,0) circle (.3) node[left]{$1$};
\draw [shading=ball, ball color=black]  (2,1) circle (.3);
\draw [shading=ball, ball color=black] (4,4) circle (.3) node[left]{$4$};
\draw [shading=ball, ball color=black]  (6,-1) circle (.3) node[right]{$5$};
\draw [shading=ball, ball color=black] (3,-3) circle (.3) node[below]{$2$};
\node at (1.7,0.3) {$3$};
\draw [line width=1.9pt] (-1,0)--(4,4);
\draw [line width=1.9pt] (4,4)--(6,-1);
\draw [line width=1.9pt] (6,-1)--(3,-3);
\draw [line width=1.9pt] (3,-3)--(-1,0);
\draw [line width=1.9pt] (3,-3)--(4,4);
\draw [line width=1pt] (-1,0)--(2,1);
\draw [line width=1pt] (2,1)--(4,4);
\draw [line width=1pt] (2,1)--(6,-1);
\draw [line width=1pt] (2,1)--(3,-3);

\end{tikzpicture}

\caption{A complex with two maximal cliques given by vertex sets $\{1,2,3,4\}$ and $\{2,3,4,5\}$.}\label{fig: BEI3cliques}
\end{figure}

\begin{example}\cite[Figure 3 (i)]{benedetti2021hamiltonian}\label{ex: DFI} Let $\Delta$ be the complex in Figure \ref{fig: BEI3cliques}. Then $J_{\Delta}$ is a unit interval DFI and corresponds to a subideal of the ideal of maximal minors of a $3\times 5$ matrix with generators $[a_1 a_2 a_3]$ indexed by facets $\{a_1,a_2,a_3\}$ of $\Delta$. Here the maximal cliques are given by the intervals $[1,4]$ and $[2,5]$.
\end{example}

\begin{notation}
Let $\Delta$ be a unit interval determinantal facet ideal and suppose that $\Delta_i = [u_i, v_i]$ and $\Delta_{j} = [u_j,v_j]$ are cliques of $\Delta$ such that $i<j$. We say $\aa \in \Delta_i\setminus \Delta_j$ if there is an integer $l$ such that $a_1,\ldots,a_l \in [u_i, v_i]\setminus( [u_i, v_i]\cap [u_j,v_j]) $. We say $\bb \in \Delta_j\setminus \Delta_i$ if there is an integer $k$ such that $b_k,\ldots,b_m \in [u_j, v_j]\setminus( [u_i, v_i]\cap [u_j,v_j]) $.
\end{notation}

Ene, Herzog, and Hibi proved that there always exists a labeling of a closed graph $G$ such that all the facets of the clique complex of $G$ correspond to intervals \cite[Theorem 2.2]{ene2011cohen}. In this way, unit interval DFIs naturally generalize binomial edge ideals of closed graphs.
In addition, unit interval DFIs have a well-understood Gr\"obner basis in a manner that further generalizes binomial edge ideals of closed graphs.

\begin{thm}\label{thm: lcmGB}\cite[Theorem 2.15]{almousa2022determinantal}
Let $\Delta$ be a pure $(m-1)$-dimensional simplicial complex on the vertex set $[n]$ admitting maximal clique decomposition $\Delta = \bigcup_{i=1}^r \Delta_i$. Let $R = \kk[X]$ be a polynomial ring over an arbitrary field $\kk$.
If the associated DFI $J_\Delta$ is unit interval, then the generators of $J_\Delta$ indexed by the facets of $\Delta$ form a reduced Gr\"obner basis with respect to \textit{any} diagonal term order, including $>$.
\end{thm}

\subsection{Rees Algebras}\label{subsec: rees}
Given a pure $(m-1)$-dimensional simplicial complex $\Delta$ and its determinantal facet ideal $J_\Delta$ in the polynomial ring $R=\mathbb{K}[x_{ij}]$, the \emph{Rees algebra} of $J_\Delta$, denoted by $\mathcal{R}(J_\Delta)$, is the graded subalgebra $R[J_\Delta \cdot t]$ of the polynomial ring $R[t]$. The \emph{Special fiber} of $J_\Delta$, denoted by $\mathcal{F}(J_\Delta)$, is the graded subalgebra $\kk[J_\Delta \cdot t]$ of the polynomial ring $\kk[t]$. Define the following standard presentations of the symmetric algebra  $\mathcal{S}(J_\Delta)$, Rees algebra $\mathcal{R}(J_\Delta)$ of $J_\Delta$, and special fiber ring $\mathcal{F}(J_\Delta)$:
\begin{align*}
\rho: R[\TT_\Delta]\longrightarrow \mathcal{S}(J_\Delta), &&
\phi: R[\TT_\Delta]\longrightarrow \mathcal{R}(J_\Delta), &&
\psi: \mathbb{K}[\TT_\Delta] &\longrightarrow \mathcal{F}(J_\Delta)
\end{align*}
where for all $i,j$, $\rho(x_{ij}) = x_{ij}=\phi(x_{ij})$, $\rho(T_\aa)=[\aa]=\psi(T_\aa)$, and $\phi(T_\aa)=[\aa]\cdot t$. 
Let $\mathcal{L}=\ker \rho$, $\mathcal{J}  = \ker \phi$, and $\mathcal{K}=\ker \psi$. The ideals $\mathcal{L}$, $\mathcal{J}$, and $\mathcal{K}$ are called the \emph{presentation ideals} of $\mathcal{S}(J_\Delta)$, $\mathcal{R}(J_\Delta)$, and $\mathcal{F}(J_\Delta)$, respectively. We sometimes refer to  the ideals $\cL$, $\cJ$, and $\cK$ as the symmetric ideal, the Rees ideal, and the special fiber ideal, respectively. When $\mathcal{L}=\mathcal{J}$, the ideal $J_\Delta$ is of \emph{linear type}. If $\mathcal{J}=\mathcal{L}+\mathcal{K}\cdot R[\TT_\Delta]$, then $J_\Delta$ is of \emph{fiber type}.

Finding the presentation ideal $\mathcal{J}$ is not easy in general. Given the presenting matrix $M$ of $J_\Delta$, the generators of $\mathcal{L}$ are given by $[\TT_{\aa^1},\dots,\TT_{\aa^{\mu}}]\cdot M$ where $\mu$ denotes the number of $(m-1)$-faces of $\Delta$. In the best scenario when $J_\Delta$ is of linear type, this gives the Rees ideal $\mathcal{J}$. However, little is known about the resolutions of determinantal facet ideals beyond the linear strand (see \cite{almousa2020resolutions}). Therefore, even finding the symmetric algebra can be difficult. The next best case is when an ideal $J_\Delta$ is of fiber type. Although the ideal of maximal minors is known to be of fiber type, this is not true in general for a determinantal facet ideal.

\begin{example}\label{ex: notFiberType}
Let $J_G$ be the binomial edge ideal corresponding to the graph $G$ with edge set $E(G)=\{(1,2),(1,4),(1,5),(2,3),(3,4),(3,5)\}$.

Then
$$f={x}_{12}{T}_{35}{T}_{14}-{x}_{14}{T}_{35}{T}_{12}-{x}_{12}{T}_{34}{T
      }_{15}+{x}_{15}{T}_{34}{T}_{12}-{x}_{14}{T}_{23}{T}_{15}+{x}_{15}{T}_{
      23}{T}_{14}$$ 
      is a minimal bihomogeneous generator of the Rees ideal, $\mathcal{J}$, of $J_G$. However, $f$ is contained in neither the symmetric ideal, $\mathcal{L}$, nor the special fiber ideal, $\mathcal{K}$, of $J_G$.
    
\end{example}

\subsection{SAGBI Bases}\label{sec: sagbi} Our goal is to use the theory of SAGBI basis deformations developed in \cite{CHV96} to find the presentation ideals of the Rees algebra and special fiber ring of unit interval determinantal facet ideals. In this way, one can use the Rees algebra of $\inn(J_\Delta)$ to understand the Rees algebra of $J_\Delta$. We recall the definition of a SAGBI basis below. For further reference on SAGBI basis theory, see \cite[Chapter 11]{sturmfels1996grobner}; for details about applications of SAGBI bases to Rees algebras, see \cite{CHV96}.

\begin{definition}\label{def:SAGBI}
Let $R$ be a polynomial ring over a field $\mathbb{K}$, and let $A\subset R$ be a finitely generated $\mathbb{K}$-subalgebra. Fix a term order $\tau$ on the monomials in $R$ and let $\inn_\tau(A)$ be the $\mathbb{K}$-subalgebra of $R$ generated by the initial monomials $\inn_\tau(a)$ where $a\in A$. We say that $\inn_\tau(A)$ is the \textit{initial algebra} of $A$ with respect to $\tau$. A set of elements $\cA\subseteq A$ is called a \emph{SAGBI basis} if $\inn_\tau(A) = \mathbb{K}[\inn_\tau(\cA)]$.
\end{definition}

\begin{definition}\label{sagbiOrder}
Let $>$ be the lexicographic order on $R = \kk[X]$ as in Notation \ref{not: initialIdeal}. Extend $>$ to a monomial order $>'$ on $R[t]$ as follows: for monomials $m_1\cdot t^i$ and $m_2\cdot t^j$ of $\mathbb{K}\left[X\right]\left[t\right]$, set $m_1\cdot t^i>' m_2\cdot t^j$ if $i>j$ or if $i=j$ and $m_1>m_2$ in $R$.
\end{definition}

The main goal of this paper is to use SAGBI basis deformation developed in \cite{CHV96} to study $\mathcal{R}(J_\Delta)$. In particular, we want to show $\inn_{>'}(\mathcal{R}(J_\Delta))=\mathcal{R}(\inn_{>}(J_\Delta))$. The first step is to understand $\mathcal{R}(\inn_{>}(J_\Delta))$. 

\begin{definition}\label{def:InitialRees}
 We define the presentations of the Rees algebra $\mathcal{R}(\inn(J_\Delta))$ and the special fiber ring $\mathcal{F}(\inn(J_\Delta))$ as follows:
\begin{align*}
\phi^\ast: R[\TT_{\Delta}] \longrightarrow \mathcal{R}(\inn(J_\Delta)), && \psi^\ast: \mathbb{K}[\TT_{\Delta}] \longrightarrow \mathcal{F}(\inn(J_\Delta))
\end{align*}
where for all $i,j$, $\phi^\ast(x_{ij}) = x_{ij}$, $\phi^\ast(T_\aa) = \xx_\aa \cdot t$, and $\psi^\ast(T_\aa)=\xx_\aa$. $\ker(\phi^\ast)$ and $\ker(\psi^\ast)$ are called \emph{presentation ideals} of $\mathcal{R}(\inn(J_\Delta))$ and $\mathcal{F}(\inn(J_\Delta))$, respectively. 
\end{definition}

In order to find a Gr\"obner basis of the defining ideal of $\mathcal{R}(\inn(J_\Delta))$, we recall a lemma that we use to achieve the goal.

\begin{lemma}\label{lem: linIndepGB}\cite[Lemma 3.1]{CHV96} Let $\kk[\bf{T}]$ be a polynomial ring equipped with a term order $\succeq$. Let $J$ be an ideal of $\kk[\bf{T}]$ and let $f_1, \dots, f_s$ be polynomials in $J$. Assume that the monomials of the set $\Omega = \{\mm \mid \mm\notin (\inn_\succeq(f_1), \dots, \inn_\succeq(f_s))\}$ are linearly independent in $\kk[{\bf{T}}]/J$. Then $f_1, \ldots, f_s$ is a Gr\"obner basis of $J$ with respect to $\succeq$.
\end{lemma}

%% file: reesClosed.tex
\section{Rees algebras of initial ideals of unit interval DFIs}\label{sec: Clique-Sorted}
In this section, we present our candidates for the defining equations of the Rees algebra of the initial ideal of a unit interval determinantal facet ideal in Definition \ref{reesEQ}. We define the \emph{clique-sorted} monomials in $R[\TT_{\Delta}]$, which will be exactly the standard monomials for $R[\TT_{\Delta}]$. We show that the candidates form a Gr\"obner basis for the presentation ideal of  $\mathcal{R}(\inn(J_\Delta))$ in Theorem \ref{thm: ReesInitialClosedDFI}.

We open this section by establishing some notation that will be assumed in the upcoming proofs.

\begin{notation}\label{notationCliques}
Let $\Delta$ be a unit interval simplicial complex on $n$ vertices with maximal clique decomposition $\Delta = \bigcup_{i=1}^r \Delta_i$ as defined in \Cref{def: unit}.
\begin{enumerate}[(i)]
 \item Let $\aa$ and $\bb$ be $m-1$-faces of $\Delta$; we say $\aa >\bb$ lexicographically if $\xx_{\aa}>\xx_{\bb}$ lexicographically as defined in Notation \ref{not: initialIdeal}.
\item Define $$\TT_{\Delta_i}=\{T_{\aa} \mid \aa \text{ is an }(m-1) \text{ face of } \Delta_i \}$$ for any $1 \leq i \leq r$.
For any monomial $\bfM \in R[\TT_{\Delta}]$, write 
$$
\bfM=\uu \cdot \mathbf{N} =\uu \cdot T_{\aa^1}\cdots T_{\aa^k}= \uu \cdot \bfG_1\cdots \bfG_r \in R[\TT_{\Delta}]
$$
where $\uu \in R$, ${\aa^1}\geq \cdots \geq {\aa^k}$ lexicographically, and $\bfG_i \in \mathbb{K}[\TT_{\Delta_i}]$.  Moreover, we impose the condition that if $T_{\aa^i}$ divides $\bfG_j$, then $T_{\aa^i} \notin \Delta _l$ for any $l<j$. This expression for $\bfM$ sorts the $T_\aa$ variables into the earliest cliques in which they appear.
\item We write $\aa \in \Delta_{i_a}$ when $\Delta_{i_a}$ is the smallest clique that contains $\aa$, i.e., $\aa \notin \Delta_l$ for any $l<i_a$.

\item  Let $\mm \in R[\inn(J_\Delta)\cdot t]$ and write $R_{\mm}=\{\bfM \in R[\TT_{\Delta}] \mid \phi^\ast(\bfM)=\mm\}$. Let $\nn \in \kk[\inn(J_{\Delta}) \cdot t]$ and write $F_{\nn}=\{\bfN \in \mathbb{K}[\TT_{\Delta}] \mid \psi^\ast(\bfN)=\nn\}$.
\end{enumerate}

\end{notation}

\begin{example}
Let $\Delta=\Delta_1 \cup \Delta_2 \cup \Delta_3=[1,6]\cup [3,7]\cup[5,9]$. Consider the monomial
$$x_{11}x_{15}x_{23}x_{28}x_{39}^2T_{124}T_{125}T_{125}T_{357}T_{367}T_{568}\in R[\TT_{\Delta}].$$ Then $\uu=x_{11}x_{15}x_{23}x_{28}x_{39}^2$, and $\bfN=T_{124}T_{125}T_{125}T_{357}T_{367}T_{568}=G_1G_2G_3$ where $G_1=T_{124}T_{125}T_{125}$, $G_2=T_{357}T_{367}$, and $G_3=T_{568}$. 
\end{example}

We define a term order for the Gr\"obner basis of the presentation ideal of $\cR(\inn(J_\Delta))$ and of  $\cF(\inn(J_\Delta))$.

\begin{definition}\label{Order}
    Let $\succ$ be the degree reverse lexicographically order on $\mathbb{K}[\TT_{\Delta}]$ and $R[\TT_{\Delta}]$ induced by $T_{\aa}>T_{\bb}$ if $\aa>\bb$ lexicographically, $T_{\aa}>x_{ij}$ for any $\aa\in \Delta$ and any $1\leq i \leq m, 1\leq j \leq n$, and $x_{ij}>x_{pq}$ if $i<p$ or $i=p$ and $j<q$. 
\end{definition}

The goal of this section is to find a Gr\"obner basis for the presentation ideal of $\cR(\inn(J_\Delta))$ and of  $\cF(\inn(J_\Delta))$ with respect to the term order $\succ$. We use a strategy similar to that found in \cite{CHV96} to achieve this goal, utilizing particular tools which are specific to our setting.

We now recall the Pl\"ucker poset from \cite[Chapter 14]{MillerSturmfels}. Endow the set
$\mathcal{P}=\{\aa\mid \aa \text{ is an $m$-subset of } [n]\}$ with a partial order as follows: if $\aa =\{a_1 <\cdots<a_m\}$ and $\bb =\{b_1 <\cdots<b_m\}$ are two elements of $\cP$, set $\aa \leq_{p} \bb$ if $a_i \leq b_i$ for all $i = 1,\dots,m$. Identifying an $m$-subset $\aa=\{a_1,\ldots,a_m\}$ with a maximal minor coming from columns $a_1,\ldots,a_m$ of an $m\times n$ matrix $X$, the poset $\mathcal{P}$ is called the \emph{Pl\"ucker poset.} 

The following lemma accomplishes two essential tasks.
First, given two facets $\aa,\bb\in \Delta$ such that $i_a\leq i_b$, we construct two other facets in $\Delta$ which we call $\min\{\aa,\bb\}$ and $\max\{\aa,\bb\}$. We show that if $\aa$ and $\bb$ are incomparable in the Pl\"ucker poset, then $\min\{\aa,\bb\}$ will be in the smaller clique, and $\max\{\aa,\bb\}$ will be in the larger clique.

\begin{lemma}\label{minMaxExist}
Adopt Notation \ref{notationCliques}. Let $\aa$ and $\bb$ be incomparable elements in the Pl\"ucker poset such that $\aa \in \Delta_{i_a}$ and $\bb \in \Delta_{i_b}$ with $i_a\leq i_b$, then
\begin{gather*}
    \cc = \{\min\{a_1,b_1\},\ldots, \min\{a_m,b_m\}\}:=\min\{\aa,\bb\}\in \Delta_{i_a}, \text{ and} \\
    \dd = \{\max\{a_1,b_1\},\ldots, \max\{a_m,b_m\}\}:=\max\{\aa,\bb\}\in \Delta_{i_b}.
\end{gather*}
\end{lemma}
\begin{proof}
 By Notation \ref{notationCliques}, it follows that $a_l \in [u_{i_a},v_{i_a}]=\Delta_{i_a}$ and $b_l \in [u_{i_b},v_{i_b}]=\Delta_{i_b}$ for all $1\leq l\leq m$. If $i_a=i_b$, then the conclusion follows trivially using the unit interval property. Suppose $i_a<i_b$, then we may assume $\aa \in \Delta_{i_a}\setminus \Delta_{i_b}$, $\bb \in \Delta_{i_b}\setminus\Delta_{i_a}$. This means we have $a_1  \in [u_{i_a},v_{i_a}]\setminus [u_{i_b},v_{i_b}]$,  $b_m  \in [u_{i_b},v_{i_b}]\setminus [u_{i_a},v_{i_a}]$, and $a_{p_{\alpha}}> b_{p_{\alpha}}$ for some $\alpha=1,...,\tau$. We notice that $a_l< b_m$ and $b_l> a_1$.
 Since $\Delta$ is a unit interval complex, and $a_{p_{\alpha}}> b_{p_{\alpha}}$, we have $d_{p_{\alpha}}=a_{p_{\alpha}} \in [u_{i_b},v_{i_b}]$ and $c_{p_{\alpha}}=b_{p_{\alpha}} \in [u_{i_a},v_{i_a}]$ for $\alpha=1,...,\tau$. 
 Therefore $\cc \in \Delta_{i_a}$ and $\dd \in \Delta_{i_b}$. 
\end{proof}

Recall that to understand the defining equations of the Rees algebra of a unit interval DFI, we first study the Rees algebra of the initial ideal of a unit interval DFI. With the lemmas above established, we are ready to present our candidates for the defining equations of Rees algebras of initial ideals of unit interval determinantal facet ideals. 

\begin{definition}\label{reesEQ}
\begin{enumerate}[(a)]
\item Let $\aa \in \Delta_{i_a}\backslash \Delta_{i_b}$ and $\bb\in \Delta_{i_b}$ be $m-1$ faces of  $\Delta$ such that $\aa>\bb$ are comparable, and $i_a<i_b$. Define polynomials of type (\ref{eq: initialKoszulClosedDFI}) to be of the form
\begin{equation}\label{eq: initialKoszulClosedDFI}
\underline{\xx_\aa T_{\bb}}-\xx_\bb T_{\aa}. 
\end{equation}
We refer to relations of type (\ref{eq: initialKoszulClosedDFI}) as \textit{Koszul-type} relations.

\item Let  $\cc = \{c_1 < c_2 < \ldots < c_{m+1}\}$ be an $m$-face of $\clique$ and $1\leq i\leq m$. Define polynomials of type (\ref{eq: initialLinSyzClosedDFI}) to be of the form \begin{equation}\label{eq: initialLinSyzClosedDFI}
\underline{x_{i c_i}T_{\cc\setminus c_i}}-x_{i c_{i+1}}T_{\cc\setminus c_{i+1}}.
\end{equation}
We refer to relations of type (\ref{eq: initialLinSyzClosedDFI}) as \textit{Eagon-Northcott-type} relations.

\item Let $\aa,\bb\in\Delta$ be incomparable elements in the Pl\"ucker poset and $\cc =\min\{\aa,\bb\}\in \Delta $ and $\dd=\max\{\aa,\bb\}\in \Delta$. 
Define polynomials of type (\ref{eq: initialPluckerClosedDFI}) to be of the form
\begin{equation}\label{eq: initialPluckerClosedDFI}
\underline{T_{\aa} T_{\bb}} - T_{\cc} T_{\dd}.
\end{equation}
We refer to relations of type (\ref{eq: initialPluckerClosedDFI}) as \textit{Pl\"ucker-type} relations.

\end{enumerate}
Let $\cG$ be the collection of binomials of (\ref{eq: initialKoszulClosedDFI}), (\ref{eq: initialLinSyzClosedDFI}), and (\ref{eq: initialPluckerClosedDFI}), and $\cG'$ be the collection of binomials of (\ref{eq: initialPluckerClosedDFI}).
\end{definition}

The following example illustrates the implicit equations listed in \Cref{reesEQ}.

\begin{example}\label{ex: BEI3cliques}
Let $\Delta$ be the unit interval simplicial complex on $6$ vertices with two maximal cliques $\Delta_1 = \{1,2,3,4,5\}$ and $\Delta_2 = \{2,3,4,5,6\}$, so $J_{\Delta}$ is generated by all three-minors $[a_1 a_2 a_3]$ of two generic $3\times 5$ matrices such that $\{a_1, a_2,a_3\} \in [1,5]$ or $\{a_1, a_2,a_3\} \in [2,6]$. We demonstrate the equations defined in Definition \ref{reesEQ} for this unit interval DFI.
\begin{enumerate}
    \item[(a)] Koszul relations between $\Delta_1$ and $\Delta_2$, e.g.,
    $$
   \underline{ \xx_{124} T_{236}} - \xx_{236}T_{124} \quad \text{ and }\quad \underline{\xx_{145} T_{256}} - \xx_{256}T_{145};
    $$
    \item[(b)] Eagon-Northcott-type relations from within $\Delta_1$ and $\Delta_2$, e.g.,
    $$
  \underline{ x_{11}T_{234}}-x_{12}T_{134}  ,\qquad \underline{x_{23}T_{256}}-x_{25}T_{236}, \; \text{ and } \;  \underline{x_{35}T_{346}}-x_{36}T_{345};
    $$
    \item[(c)] Pl\"ucker relations, e.g., 
   $$ \underline{T_{125}T_{134}}-T_{124}T_{135},   \qquad                             
 \underline{T_{145}T_{236}}-T_{135}T_{246}, \;  \text{ and } \;                              
 \underline{T_{256}T_{345}}-T_{245}T_{356}.   $$
 Notice that the middle relation arises from $\aa = \{1,4,5\}$ and $\bb =\{ 2,3,6\}$ coming from distinct maximal cliques of $\Delta$.
\end{enumerate} 
\end{example}

\begin{remark}\label{noCycle}
    One can check easily that the marked monomials of (\ref{eq: initialKoszulClosedDFI}), (\ref{eq: initialLinSyzClosedDFI}), and (\ref{eq: initialPluckerClosedDFI}) are the leading monomials with respect to $\succ$ defined in \Cref{Order}.
\end{remark}

We are ready to define \emph{clique-sorted} monomials in $R[\TT_{\Delta}]$. The clique-sorted monomials will be exactly those \emph{not} contained in the ideal generated by the underlined monomials in the polynomials of types (\ref{eq: initialKoszulClosedDFI}), (\ref{eq: initialLinSyzClosedDFI}), and (\ref{eq: initialPluckerClosedDFI}). 

\begin{notation}\label{not: addElem}
Let $\aa$ be a face of $\Delta$, and let $q\in [n]$ such that $a_{p-1}<q<a_{p}$, $p\in [m]$. Define $(\aa\cup q)$ to be the ordered tuple
$$
(a_1, \ldots, a_{p-1}, q, a_{p}, \ldots, a_m).
$$
When we write $(\aa\cup q)\backslash a_p \in \Delta$, we mean $\{a_1, \ldots, a_{p-1}, q, a_{p+1}, \ldots, a_m\} \in \Delta$.
\end{notation}

\begin{definition}\label{def: cliqueSorted} 
Adopt Notation \ref{notationCliques} and \ref{not: addElem}. A monomial $$\bfM=\uu \cdot \bfN =\uu \cdot T_{\aa^1}\cdots T_{\aa^k}=\uu \cdot G_1\cdots G_r \in R[\TT_{\Delta}]$$ is \emph{clique-sorted} if it satisfies the following properties:
\begin{enumerate}[(a)]
    \item if for all $1\leq j \leq m$ and all $1\leq i<k$ we have $a_j^i \leq a_j^{i+1}$;
    \item if $T_\aa$ divides some $G_l$, and $x_{pq}$ divides $\uu$, then either 
    $ (\aa\cup q)_p \neq q $ or $(\aa \cup q)\backslash a_p \notin \Delta$;
    \item for any $\xx_\aa$ and $T_\bb$ dividing $\bfM$ such that $\aa\in \Delta_{i_a}$ and $\bb \in \Delta_{i_b}$, we have that if $i_a \neq i_b$ and $\aa,\bb$ are comparable elements in the Pl\"ucker poset, then $i_b\leq i_a$.

\end{enumerate}
A monomial $\bfN \in \kk[\TT_\Delta]$ (so $\uu\in\kk$) satisfying property (a) is \emph{clique-sorted} in $\kk[\TT_\Delta]$.

\end{definition}

\begin{remark}\label{cliqueSortNotLead}
   If $\bfM \in R[\TT_{\Delta}]$ is clique-sorted, then none of the marked monomials of $\cG$ divide $M$. Similarly for $\bfN \in \kk[\TT_{\Delta}]$ is clique-sorted, then none of the marked monomials of $\cG'$ divide $N$. Moreover, if  $\bfN=\prod_{l=1}^r T_{\aa^l} \in \kk[\TT_{\Delta}]$ is clique-sorted then $\aa^i$ and $\aa^j$ are comparable for all $i,j$ by the definition. 
\end{remark}

We first show the existence of clique-sorted monomials.

\begin{lemma}\label{cliqueSortedExist}
Let  $\mm \in R[\inn (J_{\Delta}) \cdot t]$ and $\nn \in  \kk[\inn (J_{\Delta})\cdot t]$. Then there exist $\bfM \in R_{\mm}$ and $\bfN \in F_{\nn}$ such that $\bfM$ and $\bfN$ are clique-sorted.
\end{lemma}

\begin{proof}

    Suppose $\bfN \in F_{\nn} \subseteq \kk[\TT_{\Delta}]$ is not clique-sorted. We will show that we can find $\bfN' \in F_{\nn}$ such that $\bfN'$ is obtained by a one step reduction from $\bfN$ using Pl\"ucker-type relations in $\cG'$ and $\bfN \succ \bfN'$ where $\succ$ is the monomial ordering from Definition \ref{Order}. Since $\succ$ is a total order, this reduction process must terminate and therefore the conclusion follows.
    
   Suppose $\bfN=\prod _{l=1}^{r}T_{\bb^{l}} \in \kk[\TT]$ is not clique-sorted, then there are $\aa$ and $\bb$ such that $T_{\aa}T_{\bb} \mid \prod _{l=1}^{r}T_{\bb^{l}}$ and $\aa$ and $\bb$ are incomparable by \Cref{cliqueSortNotLead}. We may use Pl\"ucker-type relation, $T_{\aa}T_{\bb}-T_{\cc}T_{\dd}$, defined in \Cref{reesEQ}. We set $\bfN'=\bfN \frac{T_{\cc}T_{\dd}}{T_{\aa}T_{\bb}}$ to obtain $\bfN \succ \bfN'$. Notice that the existence of $\cc, \dd \in \Delta$ is proven in \Cref{minMaxExist}.

 We now consider $\bfM=\uu \cdot \prod _{l=1}^{r}T_{\aa^{l}}\in R_{\mm}$ and assume $\prod _{l=1}^{r}T_{\aa^{l}}$ is clique-sorted by the case above. Suppose $\bfM$ is not clique-sorted.  Then we have one of the possible situations: (a) there are $\aa \in \Delta$ such that $T_{\aa}\mid \prod _{l=1}^{r}T_{\aa^{l}}$ and $x_{pq} \mid \uu$ satisfying $(\aa \cup q)_p=q$ and $(\aa\cup q)\setminus a_p \in \Delta$; (b) there are $\aa,\bb\in\Delta$ where $\xx_{\aa}$ divides $\uu$ and $T_{\bb}$ divides $\prod _{l=1}^{r}T_{\aa^{l}}$ such that $\aa \in \Delta_{i_a}\backslash \Delta_{ib}$ and $\bb \in \Delta_{i_b}$ but $i_a < i_b$. 

When the condition (a)  appears, we can use Eagon-Northcott-type relation, $x_{pq}T_{\aa}-x_{pa_p}T_{(\aa\cup q)\setminus a_p}$ to obtain $\bfM \frac{x_{pa_p}T_{(\aa\cup q)\setminus a_p}}{x_{pq}T_{\aa}}$ and set $\bfM'=\bfM \frac{x_{pa_p}T_{(\aa\cup q)\setminus a_p}}{x_{pq}T_{\aa}}$, then we have $\bfM \succ \bfM'$. When the condition (b) occurs, we use Koszul-type relation, $\xx_{\aa}T_{\bb}-\xx_{\bb}T_{\aa}$ to obtain  $\bfM \frac{\xx_{\bb}T_{\aa}}{\xx_{\aa}T_{\bb}}$ and set $\bfM'=\bfM \frac{\xx_{\bb}T_{\aa}}{\xx_{\aa}T_{\bb}}$, then we have $\bfM \succ \bfM'$. Again, $\succ$ is a total order, therefore the reduction process stops. This concludes the proof of the lemma. 
\end{proof}

Before we show the uniqueness of the clique-sorted monomials in \Cref{Unique-CliqueSorted}, we need the following lemma first.

\begin{lemma}\label{SameDegree}
Adopt Notation \ref{notationCliques}. Let $\mm \in R[\inn (J_{\Delta}) \cdot t]$. Suppose $\bfM=\uu \cdot \bfG_1\cdots \bfG_r \in R_{\mm}$ and $\bfM'=\uu' \cdot \bfG'_1\cdots \bfG'_r  \in R_{\mm} $ are clique-sorted. Then $\deg\bfG_i=\deg\bfG'_i$ for all $i\in [r]$.
\end{lemma}

\begin{proof}
    For each $i$, let $\bfG_i\in F_{\gg_i}$ and $\bfG'_i\in F_{\gg'_i}$. We first notice that $\deg \bfG_1\cdots \bfG_r=\deg \bfG'_1\cdots \bfG'_r=k$ where $k$ is the degree of $t$ in $\mm$. We show $\deg\bfG_i=\deg\bfG'_i$ inductively with respect to $k$ and $r$. If $r=1$, then nothing is to be proven. We let $N_1:=\{j \mid \prod _{i=1}^j\xx_{\aa^{i}}|\mm, \aa^i \in \Delta_1 \text{ for all }i \text{, and } j\leq k\}.$ Without loss of generality, we may assume $N_1 \neq \emptyset $, and set $k_1:=\max N_1$.  We claim $k_1=\deg \bfG_1$ and similarly, $k_1=\deg \bfG'_1$. 
    We consider $\mm_1=\frac{\mm}{\gg_1\cdot t^{k_1}}$ and $\mm_1'=\frac{\mm}{\gg_1'\cdot t^{k_1}}$, then $\frac{M}{\bfG_1}=\uu\cdot \bfG_2\cdots \bfG_r \in R_{\mm_1} $ and $\frac{M}{\bfG'_1}=\uu'\cdot \bfG'_2\cdots \bfG'_r \in R_{\mm'_1}$, i.e. $\mm_1, \mm'_1 \in R[\inn(J_\Delta)\cdot t]$ and degree of $t$ in $\mm_1$ and $ \mm'_1$ are $k-k_1$. Then by induction, we have  $\deg\bfG_i=\deg\bfG'_i$ for all $i\in [r]$.

We now prove the claim. It is clear that $k_1\geq \deg \bfG_1$. Suppose $k_1>\deg \bfG_1$. Then we can find some $\xx_{\aa}$ which divides $\uu$ such that $\aa \in \Delta_{1}$. Since $r\geq 2$, we can find $T_{\bb} \mid \bfG_2\cdots \bfG_r$ such that $\bb \in  \Delta_i \setminus \Delta_1$ for some $i>1$. If $\aa \in \Delta_1 \setminus \Delta_{i}$ then we use the Koszul-type relation $\xx_{\aa}T_{\bb}-\xx_{\bb}T_{\aa}$ to reduce $\bfM$, a contradiction to $\bfM$ being clique-sorted. If $\aa \in \Delta_1 \cap \Delta_{i}$, then there is a minimal $p$ such that $a_p<b_p$. Notice that such $p$ must exist, since $\bb\in \Delta_i\setminus \Delta_1$ implies we have at least $b_m \in \Delta_i \setminus \Delta_1$ and $a_m<b_m$. If $p=1$, we have $(\bb \cup a_1)\setminus b_1 \in \Delta_i$. If $p>1$, then we have $b_{p-1}\leq a_{p-1}<a_p<b_p<b_{p+1}$, hence   $(\bb \cup a_p)\setminus b_p \in \Delta_i$ as well. Therefore we use the Eagon-Northcott-type relation $x_{pa_p}T_{\bb}-x_{pb_p}T_{(\bb\cup a_p)\setminus b_p}$ to reduce $\bfM$, which is again a contradiction to $\bfM$ being clique-sorted. This concludes the proof of the lemma. 
\end{proof}

\begin{lemma}\label{Unique-CliqueSorted}
   Adopt Notation \ref{notationCliques}. and let $\mm \in R[\inn(J_\Delta)\cdot t]$. Then there exists a unique clique-sorted monomial $\bfM \in R_{\mm}$. Similarly, let $\nn \in \kk[\inn(J_{\Delta}) \cdot t]$; then there exists a unique clique-sorted monomial $\bfN \in F_{\nn}$.
\end{lemma}

\begin{proof} 
We will prove the first statement and the second statement will follow.
By \Cref{cliqueSortedExist}, there exists $\bfM= \uu \cdot \bfG_1\cdots \bfG_r \in R_{\mm} \subseteq R[\TT_{\Delta}]$ such that $\bfM$ is clique-sorted. We give a construction to find a clique-sorted monomial $\bfM' =\uu' \cdot G_1'\cdots G_r' \in R_m$ and show that we must have $\bfM=\bfM'$. Let $k$ be the degree of $t$ in $\mm$. We show this by induction on $k$ and $r$.

Let $N_1:=\{j \mid \prod _{i=1}^j\xx_{\aa^{i}}|\mm, \aa^i \in \Delta_1 \text{ for all }i \text{, and } j\leq k\}.$ Without loss of generality, we may assume $N_1\neq \emptyset$ and set $k_1:=\max N_1$. By \Cref{SameDegree}, we have $\deg \bfG_1=k_1$. We proceed with the following steps (a), (b), and (c) to define $\bfG'_1$, and we claim $\bfG'_1=\bfG_1$.
 
    \begin{enumerate}[(a)]
    \item  Sort the $x_{1q}$ variables dividing $\mm$ by $q$ so that we have $$
x_{1q_{1}} \leq x_{1 q_{2}} \leq \cdots x_{1 q_{d_1}} \leq x_{1 v_1-(m-1)} \leq x_{1 q_{d_1+1}} \cdots \leq x_{1 q_{e_1}}.
 $$

Then set $a^{1,\tau}_1 = q_{\tau}$ for $1\leq \tau \leq k_1$. Notice $q_{\tau}$ exists for all $1\leq \tau \leq k_1$ by the definition of $k_1$.

\item Sort the $x_{2q}$ variables in the same way:
\[
x_{2 q_{1}} \leq \cdots \leq x_{2 q_{d_2}}\leq x_{2 v_1-(m-2)} \leq x_{2 q_{d_2+1}} \leq \cdots \leq x_{2 q_{e_2}}.
\]
Notice those $q_i's$ are elements inside $[n]$ and they are different from those $q_i's$ for $x_{1q}$.
Let $A^{1,\tau}_2:=\{w \mid   a^{1,\tau}_1<q_{w},  a^{1,\tau-1}_2=q_{j} \leq q_{w} \leq v_1-(m-2), \text{ and } j<w \}$. Then $A^{1,\tau}_2$ is not empty by the definition of $k_1$ again. Set $a^{1,\tau}_2 =q_{ \min A^{1,\tau}_2}$ for $1\leq \tau \leq k_1$.

\item Repeat this process, consecutively sorting all the $x_{pq}$ variables dividing  $\mm$ for a fixed $p$ and setting
$a^{1,\tau}_p = q_{\min A^{1,\tau}_p}$ where
\[
A^{1,\tau}_p:=\{w \mid  a^{1,\tau}_{p-1} < q_{w}, a^{1,\tau-1}_{p}=q_j \leq  q_{w} \leq v_1-(m-p) \text{ and } j<w  \},
\]
\end{enumerate}

Set $\gg'_1=\prod _{\tau=1}^{k_1}\xx_{\aa^{1,\tau}}$.  Notice that $\xx_{\aa} \nmid \frac{\mm}{\gg'_1}$ for any $\aa \in \Delta _1$. Let $\bfG'_1=\prod _{\tau=1}^{k_1}T_{\aa^{1,\tau}}$. We now show that $\bfG'_1=\bfG_1$. We prove this by contradiction.
 From construction, we have $\aa^{1,1}\geq \aa^{1,2} \geq \cdots \geq \aa^{1,k_1}$ lexicographically, and $\aa^{1,i}$ and $\aa^{1,j}$ are comparable for any $i,j$. We write $\bfG_s=\prod _{\tau=1}^{k_s}T_{\bb^{s,\tau}}$ for $s\in [r]$ and assume $\bb^{i,l}\geq \bb^{j,k}$ lexicographically when $i<j$ or $i=j$ and $l<k$, and $\bb^{i,l}$ and $\bb^{j,k}$ are comparable for any $i,j,l,k$ by the clique-sorted assumption of $\bfM$ and \Cref{cliqueSortNotLead}.

Suppose $\bfG'_1 \neq \bfG_1$.
Let $\beta=\min\{\tau \mid \bb^{1,\tau} \neq \aa^{1,\tau}\}$. As $\aa^{1,\beta}\neq \bb^{1,\beta}$, we must have a $p$ such that $a_{p}^{1,\beta}< b_{p}^{1,\beta}$ by the construction of $a_{p}^{1,\beta}$. Let $p$ be the minimal number such that $a_{p}^{1,\beta}< b_{p}^{1,\beta}$, i.e., $a_{p-1}^{1,\beta}= b_{p-1}^{1,\beta}$. 

Since we have comparable among $\bb^{l,\tau}$ and we order them lexicographically, we obtain $a_{p}^{1,\beta}\neq b_{p}^{l,\tau}$ for any $l>1$ or $l=1$ and $\tau>\beta$. This means $x_{{p}a_{p}^{1,\beta}}|\uu$. Moreover, we have $(\bb^{1,\beta}\cup a_{p}^{1,\beta})\backslash b_{p}^{1,\beta} \in \Delta_{1}$ and $(\bb^{1,\beta}\cup a_{p}^{1,\beta})_{p}=a_{p}^{1,\beta}$. Now we can use Eagon-Northcott-type relation, $x_{{p}a_{p}^{1,\beta}}T_{\bb^{1,\beta}}-x_{{p}b_{p}^{1,\beta}}T_{(\bb^{1,\beta}\cup a_{p}^{1,\beta})\backslash b_{p}^{1,\beta}}$, to reduce $\bfM$, a contradiction to $\bfM$ is clique-sorted. Therefore we must have $\bfG'_1 = \bfG_1$. 

If $r=1$, then $\bfM$ is the unique clique-sorted monomial with respect to $\mm$.  If $r>1$, then $\mm_1=\frac{\mm}{\gg'_1 \cdot t^{k_1}}=\frac{\mm}{\gg_1 \cdot t^{k_1}} \in R[\inn(J_{\Delta'})\cdot t]$ as we have $\bfM_1=\bfM \frac{1}{\bfG_1}\in R_{\mm_1} \subseteq R[T_{\Delta'}]$ where $\Delta'=\cup _{i=2}^r \Delta_i$. By induction on $k$ and on $r$, $\bfM_1$ is the unique clique-sorted monomial with respect to $\mm_1$ and it can be constructed using exact same steps of the construction of $\bfG_1$ repeatedly. Then $\bfM=\bfG_1 \bfM_1$ is the unique clique-sorted monomial with respect to $\mm$ and this concludes the proof of the lemma.
\end{proof}

We are now ready for the main theorem of the section.

\begin{thm}\label{thm: ReesInitialClosedDFI} 
Adopt Definition \ref{reesEQ}.
 Under the monomial order $\succ$ on the ring $R[\TT_{\Delta}]$ and $\mathbb{K}[\TT_{\Delta}]$, $\cG$ is a Gr\"obner basis of the presentation ideal of $\mathcal{R}(\inn(J_\Delta))$, and $\cG'$ is a Gr\"obner basis of the presentation ideal of  $\mathcal{F}(\inn(J_\Delta))$.
 In particular, $\inn(J_\Delta)$ is of fiber type.
\end{thm}
\begin{proof}
We apply the same strategy as in the proof of \cite[Proposition 3.2]{CHV96}. Clearly, polynomials in $\cG$ sit inside $\ker\phi^\ast$, and Pl\"ucker type relations in $\cG'$ sit inside $\ker \psi^\ast$. The term order $\succ$ on $R[\TT_{\Delta}]$ selects the underlined monomials as leading terms, and its restriction on $\kk[\TT]$ selects the underlined monomials of Pl\"ucker type relations in $\cG'$ as leading terms. Let $L$ be the ideal generated by the underlined monomials. To show that $\cG$ forms a Gr\"obner basis for $\cR(\inn(J_\Delta))$, and that Pl\"ucker type relations $\cG'$ forms a Gr\"obner basis for $\cF(\inn(J_\Delta))$, it suffices to check that all monomials not contained in $L$ are linearly independent modulo $\ker\phi^\ast$ by Lemma \ref{lem: linIndepGB}. Identify $R[\TT_{\Delta}]/\ker\phi^\ast$ with $R[\inn(J_\Delta)\cdot t]$ via the natural isomorphism induced by $\phi^\ast$. The proof is complete if we show that every monomial in $R[\inn(J_\Delta)\cdot t]$ corresponds uniquely to a clique-sorted monomial in $R[\TT_{\Delta}]$ modulo $\ker\phi^\ast$, and this is shown in \Cref{Unique-CliqueSorted}.  This completes the proof of the theorem.
\end{proof}

\section{Rees Algebras of unit interval Determinantal Facet Ideals}\label{sec: reesClosed}

The goal of this section is to find a Gr\"obner basis for the defining ideal of the Rees algebra of a unit interval determinantal facet ideal using the  Gr\"obner basis for the defining ideal of the Rees algebra of the initial ideal that we describe in Theorem \ref{thm: ReesInitialClosedDFI}.  

The motivation for the following technical lemma can be found in \cite[Section 7.1 and 7.2]{brunsherzog}. 
Building on the assumptions of \Cref{minMaxExist}, we show that a certain class of minors which we denote by $[\ee]$'s and $[\gg]$'s are guaranteed to be minimal generators of $J_\Delta$. These minors will appear in the upcoming Pl\"ucker relations in Theorem \ref{thm: sagbiLifts} and \Cref{reesThirdEq}.
  
\begin{lemma}\label{PluckersExist}
Adopt Notation \ref{notationCliques}. Let $\aa$ and $\bb$ be incomparable elements in the Pl\"ucker poset such that $\aa \in \Delta_{i_a}$ and $\bb \in \Delta_{i_b}$ with $i_a\leq i_b$. Write $\cc =\min\{\aa,\bb\}\in \Delta_{i_a}$ and $\dd=\max\{\aa,\bb\}\in \Delta_{i_b}$. Let $\ee=\{e_1,e_2,\ldots,e_m\}$ and $\gg=\{g_1,g_2,\ldots,g_m\}$ satisfy the following conditions: 
\begin{enumerate}[(a)]
    \item  $\ee \leq_p \gg$ and $\ee \leq_p \cc$ in the Pl\"ucker poset.
    \item  The sequence $(e_1,e_2,\ldots,e_m,g_1,g_2,\ldots,g_m)$ arises from the sequence\\ $(a_1,a_2,\ldots,,a_m,b_1,b_2,\ldots,b_m)$ by a permutation.
\end{enumerate}
Then $\ee \in \Delta_{i_a}$ and $\gg\in \Delta_{i_b}$.
\end{lemma}

\begin{proof}
Write $\Delta_{i_a}=[u_{i_a},v_{i_a}]$ and $\Delta_{i_b}=[u_{i_b},v_{i_b}]$. 
If $i_a=i_b$, then $a_{p},b_p\in [u_{i_a},v_{i_a}]$ for all $1\leq p\leq m$, hence $\ee, \gg \in \Delta_{i_a}$ by the requirement (b).
Therefore, we may assume $i_a<i_b$, so $u_{i_a}<u_{i_b}$ and $v_{i_a}<v_{i_b}$. 
Suppose that $a_{p}\in [u_{i_a},v_{i_a}]\setminus[u_{i_b},v_{i_b}]$ when $1\leq p\leq \tau$ and $a_{p}\in [u_{i_a},v_{i_a}]\cap [u_{i_b},v_{i_b}]$ when $\tau<p\leq m$.

This implies that $c_p=\min\{a_{p},b_p\}=a_{p}$ for all $1\leq p\leq \tau$. 
Notice that $\tau$ is the number of elements of the set $\{a_1,a_2,\ldots,a_m,b_1,b_2,\dots,b_m\}$ which are contained in $[u_{i_a},v_{i_a}]\setminus [u_{i_b},v_{i_b}]$. 
Since $\ee \leq_p \cc$, we have $e_p=a_{p} \in [u_{i_a},v_{i_a}]\setminus [u_{i_b},v_{i_b}]$ for $1\leq p \leq \tau$ and  $e_p \in [u_{i_a},v_{i_a}]\cap [u_{i_b},v_{i_b}]$ for all $\tau<p\leq m$. 
This shows $\ee \in \Delta_{i_a}$.
Moreover, $$\{a_1,a_2,\ldots,a_m,b_1,b_2,\dots,b_m\}\setminus\{e_1,\dots,e_m\} \subseteq [u_{i_b},v_{i_b}]$$
implying that $\gg\in \Delta_{i_b}$.
\end{proof}

\begin{example}\label{ExamplePluckers}
Let $\Delta$ be the unit interval simplicial complex on $6$ vertices with two maximal cliques $\Delta_1 = \{1,2,3,4,5\}$ and $\Delta_2 = \{2,3,4,5,6\}$, so $J_{\Delta}$ is generated by all three-minors $[a_1 a_2 a_3]$ of two generic $3\times 5$ matrix such that $\{a_1, a_2,a_3\} \in [1,5]$ or $\{a_1, a_2,a_3\} \in [2,6]$. 
Let $\aa=\{1,4,5\}$ and $\bb=\{2,3,6\}$. Then $\cc=\{1,3,5\}$ and $\dd=\{2,4,6\}$. Here are all possible pairs of $\ee$ and $\gg$ that are not equal to $\cc$ and $\dd$: $\ee_1=\{1,2,3\}$ and $\gg_1=\{4,5,6\}$;  $\ee_2=\{1,2,4\}$ and $\gg_2=\{3,5,6\}$; $\ee_3=\{1,2,5\}$ and $\gg_3=\{3,4,6\}$; $\ee_4=\{1,3,4\}$ and $\gg_4=\{2,5,6\}$.

\end{example}

Let $\aa=\{a_1,\ldots,a_m\}$ and $\bb=\{b_1,\ldots, b_m\}$ be $(m-1)$-dimensional faces of  $\Delta$, and $\cc=\min\{\aa,\bb\}$ and $\dd=\max\{\aa,\bb\}$ as defined in \Cref{minMaxExist}. Then $[\aa][\bb]$ can be written as a linear combination of products of minors $[\ee][\gg]$ satisfying assumptions in \Cref{PluckersExist}. (see for example the proof of Theorem 6.46 in \cite{EH2012grobner}).

We are now ready to give a SAGBI basis of the Rees algebra for a unit interval determinantal facet ideal.

\begin{thm}\label{thm: sagbiLifts}  Let $\Delta$ be a pure and unit interval $(m-1)$-dimensional simplicial complex. The polynomials of the set $\{x_{ij}\}\cup \{\left[\aa\right]\cdot t\mid \aa$ is a facet of $\Delta \}$ form a SAGBI basis of the Rees algebra $\mathcal{R}(J_\Delta)$ with respect to the monomial order $>'$ defined in Definition \ref{sagbiOrder}. In particular,
$$
\inn_{>'}(\mathcal{R}(J_\Delta)) = \mathbb{K}\left[X\right]{[}\inn_{>}(J_\Delta)\cdot t{]} = \mathcal{R}(\inn_{>}(J_\Delta)).
$$
Additionally, the polynomials of the set $\{\left[\aa \right]\mid \aa$ is a facet of $\Delta \}$ form a SAGBI basis of the $\mathbb{K}$-algebra $\mathbb{K}{[}J_\Delta{]}$ with respect to the lexicographic monomial order $>$ as in Notation \ref{not: initialIdeal}; in particular, $\inn_{>}(\mathbb{K}{[}J_\Delta{]}) = \mathbb{K}{[}\inn_{>}(J_\Delta){]}$. 
\end{thm}

\begin{proof}
Polynomials of types (\ref{eq: initialKoszulClosedDFI}), (\ref{eq: initialLinSyzClosedDFI}), and (\ref{eq: initialPluckerClosedDFI}) form a Gr\"obner basis, and therefore a (not necessarily minimal) generating set, of $\ker \phi^{\ast}$ by Theorem \ref{thm: ReesInitialClosedDFI}. It suffices to show that for any $f=f_1-f_2$ in the generating set of $\ker \phi^{\ast}$, $\phi(f)$ is a linear combination of elements of the form $\lambda \uu([\aa{]}\cdot t)^k$ with $\lambda \in \mathbb{K}\setminus \{0\}$, $k \in \mathbb{N}$, $\uu$ a monomial in the $x_{ij}$, and $\inn_{>'}(\phi (f_1)) > \inn_{>'}(\uu (\left[\mathbf{a}\right]\cdot t)^k))$; see, for example, \cite[Theorem 6.43]{EH2012grobner} or \cite[Theorem 3.3]{CHV96}.

Observe the following elementary facts: $\sum_{\substack{\p\in \mathfrak{S}_m}}\sgn(\p) x_{1\p(a_1)}\cdots x_{m\p(a_m)} =[\aa]$;\\ $\sum_{\substack{\p\in \mathfrak{S}_m}}\sgn(\p) x_{1\p(b_1)}\cdots x_{m\p(b_m)} =[\bb]$; and $[\aa][\bb]-[\bb][\aa]=0$.
Then the linear relation (\ref{eq: initialKoszulClosedDFI}) can be lifted to \begin{align*}
&\xx_\aa [\bb]\cdot t - \xx_\bb [\aa]\cdot t=\\
&\left(\sum_{\substack{\p\in \mathfrak{S}_m \\ \p\neq id}}\sgn(\p) x_{1\p(a_1)} \cdots x_{m\p(a_m)}\right)[\bb]\cdot t - \left(\sum_{\substack{\p\in \mathfrak{S}_m \\ \p\neq id}}\sgn(\p) x_{1\p(b_1)} \cdots x_{m\p(b_m)}\right)[\aa]\cdot t
\end{align*}
where $\mathfrak{S}_m$ denotes the symmetric group on $m$ letters. Observe that every monomial in $\sum\limits_{\substack{\p\in \mathfrak{S}_m \\ \p\neq id}} x_{1\p(a_1)} \cdots x_{m\p(a_m)}$ is less than $\xx_\aa$ with respect to $>$.

The linear relation (\ref{eq: initialLinSyzClosedDFI}) can be lifted to \begin{align*}
(-1)^{i+i}x_{i c_i}[\cc\setminus c_i]\cdot t +(-1)^{i+i+1} x_{i c_{i+1}}[\cc\setminus c_{i+1}]\cdot t = \sum_{\substack{j\in\{1,\ldots, m+1\} \\ j\neq i, i+1}} (-1)^{i+j+1} x_{ic_j} [\cc \setminus c_j]\cdot t
\end{align*}
where $\cc = \{c_1 < c_2 < \ldots < c_{m+1}\}$ is an $m$-face of $\clique$. This is because we have $\sum_{j\in\{1,\ldots, m+1\} } (-1)^{i+j} x_{ic_j} [\cc \setminus c_j]=0$.  
If $j<i$, then the lead monomial of $x_{i c_j} [\cc\setminus c_j]$ on the right-hand side of the equation is $$x_{ic_j}x_{1c_1}\cdots x_{j-1 c_{j-1}}x_{j c_{j+1}}\cdots x_{i c_{i}}\cdots x_{m c_{m+1}}.$$
On the other hand, the lead monomial of the left-hand side is $$x_{ic_i}x_{1c_1}\cdots x_{j-1 c_{j-1}}x_{j c_{j}}\cdots x_{i-1 c_{i-1}} x_{i c_{i+1}}\cdots x_{m c_{m+1}}.$$ Since $x_{j c_{j+1}} < x_{j c_j}$, the desired condition on the initial monomials is satisfied.

If $j>i+1$, then the leading monomial of  $x_{i c_j} [\cc\setminus c_j]$ is $$x_{ic_j}x_{1c_1}\cdots x_{i-1 c_{i-1}}x_{i c_{i}} x_{i+1 c_{i+1}}\cdots x_{j-1 c_{j-1}}x_{j c_{j+1}}\cdots x_{m c_{m+1}}.$$
The lead monomial of the left-hand side is $$x_{ic_i}x_{1c_1}\cdots x_{i-1 c_{i-1}} x_{i c_{i+1}}\cdots x_{j-1 c_{j-1}}x_{j c_{j}}\cdots x_{m c_{m+1}}.$$ The conclusion follows by $x_{i c_j} < x_{i c_{i+1}}$.

The Pl\"ucker relation (\ref{eq: initialPluckerClosedDFI}) can be lifted to the standard Pl\"ucker relation
 \begin{equation}\label{prePlucker}
[\aa][\bb]\cdot t^2 - [\cc][\dd]\cdot t^2  = \sum_{\substack{c_{\ee,\gg}\neq 0 \\ [\ee]\leq_p [\cc], [\ee]\leq_p [\gg]\neq [\dd]}} c_{\ee,\gg} \cdot [\ee] \cdot [\gg] \cdot t^2
\end{equation}
where $0\neq c_{\ee,\gg}\in \mathbb{K}$ and $[\ee], [\gg]$ are in the Pl\"ucker poset such that they satisfy the assumptions in \Cref{PluckersExist} with $\ee, \gg \in \Delta$. It is well-known that $\inn_>([\ee][\gg]) < \inn_>([\aa][\bb])$; see, for example, \cite[Theorem 6.46]{EH2012grobner}.
\end{proof}

We present the defining equations of the defining ideal of the Rees algebra of the determinantal ideal.

\begin{definition}\label{reesThirdEq} Adopt Notation \ref{notationCliques} and Definition \ref{reesEQ}. Let $\Delta$ be a unit interval and pure $(m-1)$-dimensional simplicial complex. 
\begin{enumerate}[(a)]
    \item Let $\aa$ and $\bb$ be $(m-1)$-dimensional faces contained in distinct cliques of $\Delta$. We define polynomials of type (\ref{eq: KoszulClosedDFI}) as
    \begin{equation}\label{eq: KoszulClosedDFI}
[\aa]\cdot T_\bb - [\bb]\cdot T_\aa
\end{equation}

\item Let $\cc = \{c_1, \ldots, c_{m+1}\}$ be an $m$-face of $\clique$. We define polynomials of type (\ref{eq: LinSyzClosedDFI}) as
\begin{equation}\label{eq: LinSyzClosedDFI}
\sum_{j\in\{1,\ldots, m+1\}} (-1)^{i+j} x_{ic_j} T_{\cc\setminus c_j}
\end{equation}
 
\item Let $\aa$ and $\bb$ be $(m-1)$-dimensional faces of $\Delta$. Let $\cc=\min\{\aa,\bb\},\dd=\max\{\aa,\bb\}$, and $\ee,\gg$ are faces of $\Delta$ satisfying the assumptions in \Cref{PluckersExist}. We write $i_\doot$ as the permutation $\p\in \mathfrak{S}_{2m}$ given by $(e_1,e_2,\ldots,e_m,g_1,g_2,\ldots,g_m)$ arises from the sequence $(a_1,a_2,\ldots,,a_m,b_1,b_2,\ldots,b_m)$. We define polynomials of type (\ref{eq: PluckerClosedDFI}) as
\begin{equation}\label{eq: PluckerClosedDFI}
T_{\aa}T_{\bb}-T_{\cc}T_{\dd}+\sum_{[\ee]\leq_p [\cc], [\ee]\leq_p [\gg]\neq [\dd]} \sgn(i_\doot)\cdot T_{\ee}T_{\gg}
\end{equation}
\end{enumerate}
\end{definition}

Now, applying \cite[Corollary 11.6]{sturmfels1996grobner}, we obtain our main result.

\begin{thm}\label{thm: reesClosedDFI} 
With respect to some monomial ordering $\omega'$ on the ring $R[\TT_{\Delta}]$, the presentation ideal of $\mathcal{R}(J_\Delta)$ has a Gr\"obner basis consisting of polynomials of types (\ref{eq: KoszulClosedDFI}), (\ref{eq: LinSyzClosedDFI}), and (\ref{eq: PluckerClosedDFI}).

In addition, the presentation ideal of $\mathcal{F}(J_\Delta)$ has a Gr\"obner basis given by polynomials of type (\ref{eq: PluckerClosedDFI}) with respect to some monomial order $\omega$ on $\kk[\TT]$. In particular, $J_\Delta$ is of fiber type.
\end{thm}

\begin{example}\label{ex:BEIfullPlucker}
Consider again the complex $\Delta$ from Example \ref{ex: BEI3cliques}. Lifting the defining equations of the Rees ideal of $\inn (J_{\Delta})$,  
we obtain the defining equations of $\mathcal{R}(J_{\Delta})$:
\begin{enumerate}
    \item[(a)] Koszul relations between $\Delta_1$ and $\Delta_2$, e.g.,
    $$
    [124] T_{236} - [236]T_{124} \text{ and } [145] T_{256} - [256]T_{145};
    $$
    \item[(b)] Eagon-Northcott-type relations from within $\Delta_1$ and $\Delta_2$, e.g.,
    $$
   x_{11}T_{234}-x_{12}T_{134}+x_{13}T_{124}-x_{14}T_{123}  \text{ and } 
   x_{24}T_{356}-x_{25}T_{346}-x_{23}T_{456}+x_{26}T_{345};
    $$
    \item[(c)] Pl\"ucker relations, e.g., 
   \begin{equation*}
\begin{aligned} &T_{125}T_{134}-T_{124}T_{135}+T_{123}T_{145} \text{ and }   \\                           
 &T_{145}T_{236}-T_{135}T_{246} +T_{123}T_{456} -T_{124}T_{356}+T_{125}T_{346}+T_{134}T_{256} .
\end{aligned}
\end{equation*}
 \end{enumerate} 
 Example \ref{ExamplePluckers} gives all possible $\ee$ and $\gg$ pairs for the last Pl\"ucker relation.

\end{example}

\begin{cor}\label{cor: propertiesReesClosed} Let $\Delta$ be a unit interval and pure $(m-1)$-dimensional simplicial complex with clique decomposition $\Delta = \bigcup_{i=1}^r \Delta_i$ and $J_\Delta$ be its corresponding determinantal facet ideal, and let $n_i$ be the size of the vertex set of each $\Delta_i$ in the clique decomposition. Then we have the following properties:
\begin{enumerate}[(i)]
    \item\label{item:koszulSF} $\mathcal{F}(J_\Delta)$ is Koszul.
    \item\label{item:koszulR} $\mathcal{R}(J_\Delta)$ is Koszul if $\Delta$ is a clique.
    \item\label{item:lintype} $J_\Delta$ is of linear type if and only if $n_i \leq  m+1$ for all $i$.
    \item\label{item:normal} $\mathcal{R}(J_\Delta)$ and $\mathcal{F}(J_\Delta)$ are normal Cohen-Macaulay domains. In particular, $\mathcal{R}(J_\Delta)$ and $\mathcal{F}(J_\Delta)$ have rational singularities if $\Char\mathbb{K} = 0$, and they are F-rational if $\Char\mathbb{K}>0$.
   
\end{enumerate}
\end{cor}

\begin{proof}
It is well-known that if the presentation ideal for an algebra has a quadratic Gr\"obner basis, then the algebra is Koszul; see, for instance, \cite[Theorem 6.7]{EH2012grobner}. This gives  (i) and (ii). To see (iii), observe that all relations of type (\ref{eq: PluckerClosedDFI}) come from faces of $\clique$ which are dimension $m+1$ or larger and $|\Delta_{i_a} \cap \Delta_{i_b}|\geq m+1$. By \cite[Proposition 13.15]{sturmfels1996grobner}, the semigroup rings, $\mathcal{R}(\inn_{>}(J_\Delta))$ and $\mathcal{F}(\inn_{>}(J_\Delta))$, are normal because their presentation ideals have square-free initial ideals by Theorem \ref{thm: ReesInitialClosedDFI}.

Applying \cite[Corollary 2.3]{CHV96} and Theorem \ref{thm: sagbiLifts}, we obtain (iv).
\end{proof}

\begin{acknowledgment*}
    The authors express sincere thanks to Yi-Huang Shen, Elisa Gorla, and the anonymous referee for their careful reading of the manuscript and their valuable suggestions. The first author was partially supported by the NSF GRFP under Grant No. DGE-1650441.
\end{acknowledgment*}